\documentclass[11pt,reqno]{article}
\usepackage{amsmath, amssymb, amsthm}

\newtheorem{theorem}{Theorem}

\newtheorem{lemma}[theorem]{Lemma}

\newtheorem*{general Gromov'}{Corollary \ref{general Gromov}$'$}

\def \proof {\noindent {\bf Proof.}\ \ }

\def \remarks {\noindent {\bf Remarks.}\ \ }

\def \endproof {{\mbox{}\nolinebreak\hfill\rule{2mm}{2mm}\par\medbreak}}

\def \R {\mathbb{R}}

\def \P {\mathbb{P}}

\def \OO {\mathcal{O}}

\def \e {\varepsilon}

\def \d {\delta}

\def \s {\sigma}

\def \< {\langle}
\def \> {\rangle}

\def \diam {{\rm diam}}

\def \rank {{\rm rank }}

\def \id {{\it id}}

\def \Lip {{\rm Lip}}

\begin{document}
\title{On random intersections of two convex bodies \\
{\normalsize Appendix to: ``Isoperimetry of waists 
and local versus global asymptotic convex geometries'' 
by R. Vershynin}}


\author{Mark Rudelson \and Roman Vershynin}
\date{}

\maketitle

In \cite{V} it is proved that the existence of nicely bounded sections
of two symmetric convex bodies $K$ and $L$ in $\R^n$ (of dimensions $k$ and $n-ck$)
implies that the random intersection $K \cap UL$ is nicely bounded
with high probability, where $U$ is a random unitary operator. 
Namely, the diameter of $K \cap UL$ is at most $C^{n/k}$ times
the larger of the diameters of the two sections, with probability 
at least $1 - e^{-n}$.

In this appendix we show how to improve the exponential bound $C^{n/k}$
to a polynomial bound, say $C(n/k)^2$. The cost for this is decreasing the 
probability from $1 - e^{-n}$ to $1 - e^{-k}$.

\begin{theorem}                     \label{two bodies ap}
  Let $0 < a < 0.03$. 
  Assume that two symmetric convex bodies $K$ and $L$ in $\R^n$
  have sections of dimensions at least $k$ and $n-ak$ respectively 
  whose diameters are bounded by $1$.
  Then for every $t > C(n/k)^{Ca}$ the random orthogonal 
  operator $U \in \OO(n)$ satisfies 
  $$
  \P \Big\{ \diam(K \cap UL) > tn/k \Big\}  <  (ct)^{-k/16}.
  $$
\end{theorem}

\remarks
{\bf 1.} Theorem 1.1 of \cite{V} is a partial case of this theorem for 
$t = C_1^{n/k}$:
$$
\P \Big\{ \diam(K \cap UL) > C^{n/k} \Big\}  <  e^{-n}.
$$

{\bf 2.} To obtain a polynomial bound on the diameter, one can choose 
$t = C_1(n/k)$ in Theorem \ref{two bodies ap} to get
$$ 
\P \Big\{ \diam(K \cap UL) > C_1(n/k)^2 \Big\}  
<  (c C_1 n/k)^{-k/16} 
< e^{-k}
$$
for an appropriate absolute constant $C_1$.
To summarize, 
\begin{quote}
  In Theorem 1.1 of \cite{V}, the body $K \cap UL$ has diameter 
  bounded by $C_1(n/k)^2$ with probability at least $1 - e^{-k}$.
\end{quote} 


A new ingredient in the proof of Theorem \ref{two bodies ap} 
is the following covering lemma.

\begin{lemma}                       \label{covering}
  Let $K$ be a convex body in $\R^n$ such that
  $K \supseteq \d D$ for some $\d > 0$. 
  Assume that there exists an orthogonal projection $P$ 
  with $\rank P = n-k$ and such that $PK \supseteq PD$.
  Then 
  $$
  N(D, 4K)  \le  (C/\d)^{2k}.
  $$
\end{lemma}

\proof
Denote the range of $P$ by $E$.
Let $f : PD \to \R^n$ be a lifting of $f$, i.e. a map such that
\begin{equation}                    \label{lifting}
f(PD) \subset K
\ \ \ \text{and} \ \ \
V := (\id - f)(PD)  \subset E^\perp.
\end{equation}
Since $V \subset PD - f(PD) \subset D - f(PD)$, 
the assumptions on $K$ and \eqref{lifting} imply that
$V \subset (\frac{1}{\d} + 1) K \cap E^\perp$. Then  
by the standard volumetric argument we have 
$$
N(V,K) \le (C/\d)^k
$$
as $\dim E^\perp = k$.
This will allow us to cover $PD$. Indeed, by \eqref{lifting}, 
$$
PD \subset f(PD) + V \subset K+V.
$$
By the submultiplicative property 
$N(K_1+K_2, D_1+D_2) \le N(K_1, D_1) \; N(K_2, D_2)$, 
which is valid for all sets $K_1, K_2, D_1, D_2$, we have 
\begin{equation}                    \label{cover 1}
N(PD, 2K) \le N(K+V, 2K) \le N(V,K) \le (C/\d)^k.
\end{equation}
Also by the assumption on $K$ and 
by the standard volumetric argument already used above,
\begin{equation}                    \label{cover 2}
N(D \cap E^\perp, K) \le N(D \cap E^\perp, \d D \cap E^\perp)
\le (C/\d)^k.
\end{equation}
Since $D \subset PD + D \cap E^\perp$, we have by the 
submultiplicative property that
$$
N(D, 3K) \le N(PD, 2K) \; N(P \cap E^\perp, K)
$$
and we finish by applying \eqref{cover 1} and \eqref{cover 2}.
\endproof

\qquad

\noindent {\bf Proof of Theorem \ref{two bodies ap}. }
We start the proof as in \cite{V}, by dualizing the statement and 
assuming that there exist orthogonal projections $P$ and $Q$ with 
$\rank P = k$ and $\rank Q = n-ak$, and such that 
\begin{equation}                        \label{PQ}
  PK \supseteq PD, \ \ \ QL \supseteq QD.
\end{equation}
Then we must prove that for $t$ as in the theorem, 
\begin{equation}                        \label{WTS}
  \P \{ (k/tn) D \subseteq K + UL \} \ge 1 - (ct)^{-k/16}.
\end{equation}
Let $\e > 0$ and let 
$$
\d_K = \sqrt{1 - \frac{\e^2 k}{n}}.
$$
By Proposition 3.1 
and Corollary 2.6 (ii) 
of \cite{V},  
$$
\s_{n-1}(K + \d_K D)
\ge  \s^\Lip_{n-1,k-1}(\sin^{-1} \d_K)\\
\ge 1 - (C \e_K)^{k/4}.
$$
Let $0 < \d_L < 1$ be a parameter.
By Lemma \ref{covering} applied to the body $L + \d_L$, 
$$
N(D, 4(L + \d_L D)) \le (C/\d_L)^{2ak}.
$$
Writing this covering number as $N(2D, 8L + 8\d_L D)$, we 
apply Lemma 4.1 
of \cite{V}. It states that if 
$\d_K + 8 \d_L < 1$ then the inclusion 
\begin{equation}                        \label{desired incl ap}
  (1 - \d_K - 8 \d_L) D  \subseteq K + 8UL
\end{equation}
holds with probability at least 
\begin{equation}                        \label{probability}
  1 - (C/\d_L)^{2ak} (C\e)^{k/4}.
\end{equation}
To finish the proof, we need to bound below 
the radius $1 - \d_K - 8 \d_L$ in \eqref{desired incl ap}
and the probability \eqref{probability}.
Since $\sqrt{1-x} \le 1 - x/2$ for $0 < x < 1$, we set
$$
\d_L = \frac{\e^2 k}{32n}
$$
to obtain 
$$
1 - \d_K - 8 \d_L  \ge  \frac{\e^2 k}{4n}.
$$
It remains to estimate the probability \eqref{probability}.
If we require that 
\begin{equation}                        \label{epsilon small}
  \e \le c_0 (k/n)^{C_0 a}
\end{equation}
for suitable absolute constants $c_0, C_0 > 0$, then 
$(C/\d_L)^{2ak} < (C\e)^{-k/8}$, hence the probability
$$
\text{\eqref{probability}} \ge 1 - (C\e)^{k/8}.
$$ 
We have thus proved that if $\e > 0$ satisfies \eqref{epsilon small}
then 
$$
\P \{ (\e^2 k/32n) D \subseteq K + UL \} 
\ge \P \{ (\e^2 k/4n) D \subseteq K + 8UL \} 
\ge 1 - (C\e)^{k/8}.
$$
It remains to set $\e^2/32 = 1/t$, and the proof is complete.
\endproof

\end{document}